\documentclass[10pt]{alggeom}
\usepackage{amsfonts}
\usepackage{mathrsfs}
\usepackage{amscd}
\usepackage{amsmath}
\usepackage{amssymb}
\usepackage{latexsym}
\usepackage{lscape}
\usepackage{comment}
\usepackage{amscd}
\usepackage{wasysym}  
\usepackage{tikz}
\usetikzlibrary{matrix,arrows}
\usepackage{tikz-cd}
\usepackage{appendix}
\usepackage[nice]{nicefrac}
\usepackage{dsfont, stmaryrd}

\usepackage[pagebackref,hyperindex,linktocpage=true]{hyperref}
\hypersetup{
    colorlinks,
    linkcolor={blue!20!black},
    citecolor={blue!20!black},
    urlcolor={blue!20!black}
}

\usepackage{color}

\usetikzlibrary{decorations.markings}

\makeatletter
\tikzcdset{
  open/.code     = {\tikzcdset{hook, circled};},
  closed/.code   = {\tikzcdset{hook, slashed};},
  open'/.code    = {\tikzcdset{hook', circled};},
  closed'/.code  = {\tikzcdset{hook', slashed};},
  circled/.code  = {\tikzcdset{markwith = {\draw (0,0) circle (.375ex);}};},
  slashed/.code  = {\tikzcdset{markwith = {\draw[-] (-.4ex,-.4ex) -- (.4ex,.4ex);}};},
  markwith/.code ={
    \pgfutil@ifundefined%
    {tikz@library@decorations.markings@loaded}%
    {\pgfutil@packageerror{tikz-cd}{You need to say %
      \string\usetikzlibrary{decorations.markings} to use arrows with markings}{}}{}%
    \pgfkeysalso{/tikz/postaction = {
      /tikz/decorate,
      /tikz/decoration={markings, mark = at position 0.5 with {#1}}}
    }
  },
}
\makeatother

\renewcommand{\geq}{\geqslant}
\renewcommand{\leq}{\leqslant}

\newtheorem{thm}{Theorem}[section]

\newtheorem{propo}[thm]{Proposition}
\newtheorem{lem}[thm]{Lemma}
\newtheorem{sublem}[thm]{Sublemma}
\newtheorem{lem-def}[thm]{Lemma-Definition}
\newtheorem{cor}[thm]{Corollary}
\newtheorem{conject}[thm]{Conjecture}
\newtheorem{propert}[thm]{Properties}
\newtheorem{observ}[thm]{Observation}

\theoremstyle{definition}
\newtheorem*{ack}{Acknowledgement}

\newtheorem{ex}[thm]{Example}

\newtheorem{rmk}[thm]{Remark}
\newtheorem{dfn}[thm]{Definition}
\newtheorem{quest}[thm]{Question}
\newtheorem{expec}[thm]{Expectation}
\newtheorem*{abs}{Abstract}

\numberwithin{equation}{section}

\newcommand{\nc}{\newcommand}

\nc{\theo}{\begin{thm}} \nc{\xtheo}{\end{thm}}
\nc{\prop}{\begin{propo}} \nc{\xprop}{\end{propo}}
\nc{\lemm}{\begin{lem}} \nc{\xlemm}{\end{lem}}
\nc{\sublemm}{\begin{sublem}} \nc{\xsublemm}{\end{sublem}}
\nc{\lemmdefi}{\begin{lem-def}} \nc{\xlemmdefi}{\end{lem-def}}
\nc{\coro}{\begin{cor}} \nc{\xcoro}{\end{cor}}
\nc{\conj}{\begin{conject}} \nc{\xconj}{\end{conject}}
\nc{\proper}{\begin{propert}} \nc{\xproper}{\end{propert}}
\nc{\obse}{\begin{observ}} \nc{\xobse}{\end{observ}}
\nc{\ques}{\begin{quest}} \nc{\xques}{\end{quest}}
\nc{\expe}{\begin{expec}} \nc{\xexpe}{\end{expec}}

\nc{\ackn}{\begin{ack}} \nc{\xackn}{\end{ack}}
\nc{\exam}{\begin{ex}} \nc{\xexam}{\end{ex}}
\nc{\rema}{\begin{rmk}} \nc{\xrema}{\end{rmk}}
\nc{\defi}{\begin{dfn}} \nc{\xdefi}{\end{dfn}}
\nc{\abst}{\begin{abs}} \nc{\xabst}{\end{abs}}

\nc{\pf}{\begin{proof}} \nc{\xpf}{\end{proof}}

\nc{\on}{\operatorname}
\nc{\fraka}{{\mathfrak a}} \nc{\bba}{{\mathbf a}}
\nc{\frakb}{{\mathfrak b}}
\nc{\frakc}{{\mathfrak c}}
\nc{\frakd}{{\mathfrak d}}
\nc{\frake}{{\mathfrak e}}
\nc{\frakf}{{\mathfrak f}}
\nc{\frakg}{{\mathfrak g}}
\nc{\frakh}{{\mathfrak h}}
\nc{\fraki}{{\mathfrak i}}
\nc{\frakj}{{\mathfrak j}}
\nc{\frakk}{{\mathfrak k}}
\nc{\frakl}{{\mathfrak l}}
\nc{\frakm}{{\mathfrak m}}
\nc{\frakn}{{\mathfrak n}}
\nc{\frako}{{\mathfrak o}}
\nc{\frakp}{{\mathfrak p}}
\nc{\frakq}{{\mathfrak q}}
\nc{\frakr}{{\mathfrak r}}
\nc{\fraks}{{\mathfrak s}}
\nc{\frakt}{{\mathfrak t}}
\nc{\fraku}{{\mathfrak u}}
\nc{\frakv}{{\mathfrak v}}
\nc{\frakw}{{\mathfrak w}}
\nc{\frakx}{{\mathfrak x}}
\nc{\fraky}{{\mathfrak y}}
\nc{\frakz}{{\mathfrak z}}
\nc{\frakA}{{\mathfrak A}}
\nc{\frakB}{{\mathfrak B}}
\nc{\frakC}{{\mathfrak C}}
\nc{\frakD}{{\mathfrak D}}
\nc{\frakE}{{\mathfrak E}}
\nc{\frakF}{{\mathfrak F}}
\nc{\frakG}{{\mathfrak G}}
\nc{\frakH}{{\mathfrak H}}
\nc{\frakI}{{\mathfrak I}}
\nc{\frakJ}{{\mathfrak J}}
\nc{\frakK}{{\mathfrak K}}
\nc{\frakL}{{\mathfrak L}}
\nc{\frakM}{{\mathfrak M}}
\nc{\frakN}{{\mathfrak N}}
\nc{\frakO}{{\mathfrak O}}
\nc{\frakP}{{\mathfrak P}}
\nc{\frakQ}{{\mathfrak Q}}
\nc{\frakR}{{\mathfrak R}}
\nc{\frakS}{{\mathfrak S}}
\nc{\frakT}{{\mathfrak T}}
\nc{\frakU}{{\mathfrak U}}
\nc{\frakV}{{\mathfrak V}}
\nc{\frakW}{{\mathfrak W}}
\nc{\frakX}{{\mathfrak X}}
\nc{\frakY}{{\mathfrak Y}}
\nc{\frakZ}{{\mathfrak Z}}
\nc{\bbA}{{\mathbb A}}
\nc{\bbB}{{\mathbb B}}
\nc{\bbC}{{\mathbb C}}
\nc{\bbD}{{\mathbb D}}
\nc{\bbE}{{\mathbb E}}
\nc{\bbF}{{\mathbb F}} \nc{\bbf}{{\mathbf f}}
\nc{\bbG}{{\mathbb G}}
\nc{\bbH}{{\mathbb H}}
\nc{\bbI}{{\mathbb I}}
\nc{\bbJ}{{\mathbb J}}
\nc{\bbK}{{\mathbb K}}
\nc{\bbL}{{\mathbb L}}
\nc{\bbM}{{\mathbb M}}
\nc{\bbN}{{\mathbb N}}
\nc{\bbO}{{\mathbb O}}
\nc{\bbP}{{\mathbb P}}
\nc{\bbQ}{{\mathbb Q}}
\nc{\bbR}{{\mathbb R}}
\nc{\bbS}{{\mathbb S}}
\nc{\bbT}{{\mathbb T}}
\nc{\bbU}{{\mathbb U}}
\nc{\bbV}{{\mathbb V}}
\nc{\bbW}{{\mathbb W}}
\nc{\bbX}{{\mathbb X}}
\nc{\bbY}{{\mathbb Y}}
\nc{\bbZ}{{\mathbb Z}}
\nc{\calA}{{\mathcal A}}
\nc{\calB}{{\mathcal B}}
\nc{\calC}{{\mathcal C}}
\nc{\calD}{{\mathcal D}}
\nc{\calE}{{\mathcal E}}
\nc{\calF}{{\mathcal F}}
\nc{\calG}{{\mathcal G}}
\nc{\calH}{{\mathcal H}}
\nc{\calI}{{\mathcal I}}
\nc{\calJ}{{\mathcal J}}
\nc{\calK}{{\mathcal K}}
\nc{\calL}{{\mathcal L}}
\nc{\calM}{{\mathcal M}}
\nc{\calN}{{\mathcal N}}
\nc{\calO}{{\mathcal O}}
\nc{\calP}{{\mathcal P}}
\nc{\calQ}{{\mathcal Q}}
\nc{\calR}{{\mathcal R}}
\nc{\calS}{{\mathcal S}}
\nc{\calT}{{\mathcal T}}
\nc{\calU}{{\mathcal U}}
\nc{\calV}{{\mathcal V}}
\nc{\calW}{{\mathcal W}}
\nc{\calX}{{\mathcal X}}
\nc{\calY}{{\mathcal Y}}
\nc{\calZ}{{\mathcal Z}}

\nc{\scrA}{{\mathscr A}}
\nc{\scrE}{{\mathscr E}}
\nc{\scrR}{{\mathscr R}}

\nc{\Bmu}{\mbox{$\raisebox{-0.59ex}{$l$}\hspace{-0.18em}\mu\hspace{-0.88em}\raisebox{-0.98ex}{\scalebox{2}{$\color{white}.$}}\hspace{-0.416em}\raisebox{+0.88ex}{$\color{white}.$}\hspace{0.46em}$}{}}

\nc{\bnu}{{\bar{ \nu}}}

\nc{\olO}{\bar{\calO}}

\nc{\al}{{\alpha}} 
\nc{\be}{{\beta}}
\nc{\ga}{{\gamma}} \nc{\Ga}{{\Gamma}}
 \nc{\hGa}{\hat{\Gamma}}
\nc{\ve}{{\varepsilon}} 
\nc{\la}{{\lambda}} \nc{\La}{{\Lambda}}
\nc{\om}{\omega} \nc{\Om}{\Omega} 
\nc{\sig}{{\sigma}} \nc{\Sig}{{\Sigma}}

\nc{\tnb}{\psi_{\rm tame}}
\nc{\oM}{\overline{{M}}}
\nc{\op}{{\on{op}}}
\nc{\ad}{{\on{ad}}}
\nc{\alg}{{\on{alg}}}
\nc{\Ad}{{\on{Ad}}}
\nc{\Adm}{{\on{Adm}}} \nc{\aff}{{\on{aff}}}
\nc{\Aut}{{\on{Aut}}}
\nc{\Bun}{{\on{Bun}}}
\nc{\cha}{{\on{char}}}
\nc{\der}{{\on{der}}}
\nc{\Der}{{\on{Der}}}
\nc{\diag}{{\on{diag}}}
\nc{\End}{{\on{End}}}
\nc{\Fl}{{\calF\!\ell}}
\nc{\Tr}{{\on{Transp}}}
\nc{\TR}{{\calT\!\calR}}
\nc{\Gal}{{\on{Gal}}}
\nc{\Gr}{{\on{Gr}}}
\nc{\rH}{{\on{H}}}
\nc{\Hom}{{\on{Hom}}}
\nc{\IC}{{\on{IC}}}
\nc{\id}{{\on{id}}}
\nc{\Id}{{\on{Id}}}
\nc{\ind}{{\on{ind}}}
\nc{\Ind}{{\on{Ind}}}
\nc{\Lie}{{\on{Lie}}}
\nc{\Pic}{{\on{Pic}}}
\nc{\pr}{{\on{pr}}}
\nc{\Res}{{\on{Res}}}
\nc{\res}{{\on{res}}} \nc{\Sat}{{\on{Sat}}}
\nc{\s}{{\on{sc}}}
\nc{\drv}{{\on{der}}}
\nc{\sgn}{{\on{sgn}}}
\nc{\Spec}{{\on{Spec}}}\nc{\Spf}{\on{Spf}} 
\nc{\Sph}{\on{Sph}}
\nc{\St}{{\on{St}}}
\nc{\tr}{{\on{tr}}}
\nc{\Mod}{{\mathrm{-Mod}}}
\nc{\Hilb}{{\on{Hilb}}} 
\nc{\Ext}{{\on{Ext}}} 
\nc{\vs}{{\on{Vec}}}
\nc{\ev}{{\on{ev}}}
\nc{\nO}{{\breve{\calO}}}
\nc{\tS}{{\tilde{S}}}
\nc{\spe}{{\on{sp}}}
\nc{\loc}{{\on{loc}}}
\nc{\Sym}{{\on{Sym}}}
\nc{\Cone}{{\on{C}}}
\nc{\syn}{{\on{syn}}}
\nc{\reg}{{\on{reg}}}
\nc{\colim}{{\on{colim}}}
\nc{\Norm}{{\on{N}}}

\nc{\nscrR}{{\mathscr{R}^{\on{nr}}}}

\nc{\GL}{{\on{GL}}}
\nc{\U}{{\on{U}}}
\nc{\Gl}{\on{Gl}} 
\nc{\GSp}{{\on{GSp}}}
\nc{\gl}{{\frakg\frakl}}
\nc{\SL}{{\on{SL}}} 
\nc{\SU}{{\on{SU}}} 
\nc{\SO}{{\on{SO}}}
\nc{\PGL}{{\on{PGL}}}

\nc{\Conv}{{\on{Conv}}}
\nc{\Rep}{{\on{Rep}}}
\nc{\Dom}{{\on{Dom}}}
\nc{\red}{{\on{red}}}
\nc{\act}{{\on{act}}}
\nc{\nr}{{\on{nr}}}
\nc{\ctf}{{\on{ctf}}}

\nc{\str}{{\on{-}}} 
\nc{\os}{{\bar{s}}}
\nc{\oeta}{{\bar{\eta}}}

\nc{\hookto}{\hookrightarrow}
\nc{\longto}{\longrightarrow}
\nc{\leftto}{\leftarrow}
\nc{\onto}{\twoheadrightarrow}
\nc{\lonto}{\twoheadleftarrow}

\nc{\uG}{{\underline{G}}}
\nc{\uA}{{\underline{A}}}
\nc{\uS}{{\underline{S}}}
\nc{\uT}{{\underline{T}}}
\nc{\uM}{{\underline{M}}}
\nc{\uP}{{\underline{P}}}
\nc{\uB}{{\underline{B}}}
\nc{\uN}{{\underline{N}}}

\nc{\ucG}{{\underline{\calG}}}
\nc{\ucA}{{\underline{\calA}}}
\nc{\ucS}{{\underline{\calS}}}
\nc{\ucT}{{\underline{\calT}}}
\nc{\ucalM}{{\underline{\calM}}}
\nc{\ucP}{{\underline{\calP}}}
\nc{\ucalN}{{\underline{\calN}}}

\nc{\bF}{{\breve{F}}}

\nc{\oFl}{{\overline{\Fl}}} 
\nc{\bU}{{\overline{U}}}
\nc{\tGr}{{\tilde{\Gr}}}
\nc{\cGr}{\calG\! r}
\nc{\oGr}{\overline{\on{Gr}}} 
\nc{\ocGr}{\overline{\calG\! r}}
\nc{\co}{{\colon}}
\nc{\sch}[1]{(Sch/{#1})}
\nc{\HypLoc}[1]{HypLoc({#1})}

\nc{\ohtimes}{\stackrel{!}{\otimes}}
\nc{\boxtilde}{\widetilde{\boxtimes}}
\nc{\vstar}{{\varhexstar}}

\nc{\Div}{\on{Div}}
\nc{\Sht}{\on{Sht}}
\nc{\Frob}{\on{Frob}}

\nc{\x}{\times}
\nc{\bsl}{\backslash}
\nc{\algQl}{{\bar{\bbQ}_\ell}}
\nc{\sF}{{\bar{F}}}
\nc{\nF}{{\breve{F}}}
\nc{\nW}{{W^{\on{nr}}}}
\nc{\sk}{{\bar{k}}}
\nc{\cont}{\on{c}}
\nc{\Supp}{\on{Supp}}
\nc{\blt}{\bullet}  
\nc{\dom}{\on{dom}}
\nc{\scon}{{\on{sc}}} 
\nc{\Affine}{\on{Aff}} 
\nc{\nscrA}{\mathscr{A}^{\on{nr}}} 
\nc{\nfraka}{{\bbf^{\on{nr}}}}
\nc{\ran}{{\rangle}}
\nc{\lan}{{\langle}}
\nc{\bk}{{\bar{k}}}
\nc{\tF}{{\tilde{F}}}
\nc{\sS}{{\bar{S}}}
\nc{\LG}{{^\text{L}\hspace{-0.04cm}G}}
\nc{\LL}{{^\text{L}\hspace{-0.07cm}L}}
\nc{\et}{{\text{\rm \'et}}}
\nc{\inv}{{\on{inv}}}
\nc{\Hecke}{{\on{Hecke}}}
\nc{\Isom}{{\on{Isom}}}
\nc{\oSht}{{\overline{\on{Sht}}}}
\nc{\umu}{{\underline \mu}}
\nc{\AIJ}{{\calO_X[{\scriptstyle{\calI\over \calJ}}]}}
\nc{\Proj}{{\on{Proj}}}
\nc{\Bl}{{\on{Bl}}}

\nc{\Pos}{{\on{Pos}}}
\nc{\Sets}{{\on{Sets}}}
\nc{\AffSch}{{\on{AffSch}}}
\nc{\Groups}{{\on{Groups}}}
\nc{\Gpds}{{\on{Groupoids}}}
\nc{\Sch}{{\on{Sch}}}
\nc{\fl}{{\on{flat}}}

\nc{\pot}[1]{ [\hspace{-0,5mm}[ {#1} ]\hspace{-0,5mm}] }
\nc{\rpot}[1]{ (\hspace{-0,7mm}( {#1} )\hspace{-0,7mm}) }

\nc{\defined}{\hspace{0.1cm}\stackrel{\text{\tiny \rm def}}{=}\hspace{0.1cm}}

\title{Algebraic Magnetism Invariants of Self-Actions of Diagonalizable Monoid Schemes}

\author{Arnaud Mayeux}

\email{arnaud.mayeux@mail.huji.ac.il}

\address{Einstein Institute of Mathematics, Edmond J. Safra Campus,
The Hebrew University of Jerusalem,
Givat Ram, Jerusalem, 9190401, Israel }


\begin{document}
\maketitle

\begin{flushleft}
\textbf{Abstract:} We provide a method to compute the pure magnets of the action of a diagonalizable monoid scheme on itself. This is described in terms of minimal generators of the sharp monoid obtained quotienting by the face of invertible elements. In particular, in this example, algebraic magnetism detects sharpness and minimal generators of the monoid modulo invertible elements.
\end{flushleft}

\section{Introduction}

\subsection{Results}
Algebraic Magnetism is a formalism providing invariants associated with an arbitrary action of a diagonalizable group (or monoid) scheme on a scheme \cite{Ma} (cf. also \cite{Ma25,Ma25'} for short introductions). In the case of the adjoint action of a maximal split torus on a reductive group, algebraic magnetism recovers the root system and some associated fundamental objects (cf. \cite{Ma25,Ma25',Ma}). In general, the theory relies heavily on commutative monoids and diagonalizable monoid schemes. Indeed if $D(Z)_S$ is a diagonalizable group scheme (attached to a base scheme $S$ and an abelian group $Z$, cf. \cite{Gr63}) acting via an action $a$ on a scheme $X$ separated over $S$, then the theory provides a decomposition $m(a) = \sqcup _{N \in \mho (a)} m^N (a)$ where $m(a)$ is the set of all submonoids of $Z$ (i.e. the so-called lattice of submonoids in semigroup theory) and $\mho (a)\subset m(a)$ is called the set of pure magnets $(N $ is the minimal element in $m^N (a)$, for inclusions of submonoids). To each submonoid $N $ of $Z$ (called a magnet), the theory provides a scheme $X^N $ together with an embedding (a monomorphism) into $X$. The scheme $X^N$ depends only on the class of $N$ in the above decomposition. 
In general $X^N$ is defined as the contravariant functor sending a scheme $T $ over $S$ to the set of equivariant morphisms $\Hom^{D(Z)_T }_{T} (A(N)_T , X_T )$ where $A(N)_T $ is the diagonalizable monoid scheme defined as $\Spec (\bbZ [N]) \times _{\bbZ } T $. So diagonalizable monoid schemes $A(N)_S$ are at the heart of the theory. 
On the other hand, it makes sense to consider the case where $X$ itself is a diagonalizable monoid scheme $A(M)_S$ for some submonoid $M$ of $Z$. At least at the internal foundational level of Algebraic Magnetism, it is very natural to understand the pure magnets in this very special case. In this paper, we answer this theoretical question. Let $M^*$ be the face of invertible elements of $M$, i.e. $M^*= \{x \in M | \exists y \in M , x+y =0 \}$. In the following $\mathscr{P}(-)$ means the set of subsets of a set (i.e. the power set), and $[-\rangle$ the submonoid generated by a subset of a given monoid. 
\begin{flushleft}
\textbf{Theorem.} Assume that $E$ is a minimal set of generators of the sharp quotient monoid $M/M^*$.
 If $M^*=0$, then $(B \subset E) \mapsto [B\rangle $ provides a bijection between  $\mathscr{P} (E)$ and $\mho (a)$. 
If $M^*\ne 0$, then $\mho (a) =\{0\} \sqcup \{ f^{-1} ([B\rangle)| B \in \mathscr{P}(E)\}$, where $f$ is the projection $ M \to M/M^*$. 
In particular $\#\mho(a)= 2^{\#E}$ if $M^*=0$ and $\#\mho(a)= 2^{\#E}+1 $ if $M^*\ne 0$.
\end{flushleft}

\subsection{Relation to other works} \label{sso}

Up to \cite{Ma25, Ma25', Ma}, this paper is mostly self-contained. However, monoids, algebraic monoids, and their actions were, of course, studied a long time ago. 

Monoids are fundamental and have been studied continuously for a long time \cite{Re40, Gre51, Ho09, How95, Og18}. In this paper, we only use commutative and cancellative monoids. In particular, the analysis of Green's relations on generators \cite{Gre51} (recall that our main result is formulated in terms of generators) simplifies dramatically. 

M. S. Putcha and L. E. Renner pioneered the study of algebraic monoids and semigroups, publishing hundreds of papers as sole authors, among the early ones being \cite{Pu80,  Pu81', Pu82,  Pu87,  Pu88, Pu93} and  \cite{Re82, Re84, Re85, Re05, Re09}. This field was then further advanced by several authors \cite{ Br14, Br14',   Hu96,  Li18, NNO20, Og18, Ri07}.
We do not aim to provide a complete historical overview here, readers can find many other recent references on algebraic semigroups elsewhere. 

Some aspects of diagonalizable monoid schemes in specific settings are used and studied in many places in the literature (toric geometry, log geometry), sometimes independently \cite{NNO20, Og18}. The reference \cite{NNO20} focuses on the equivalence between the opposite category of commutative monoids and that of commutative monoid $k$-schemes that are diagonalizable, for any field $k$ (compare with \cite{Gr63} for commutative groups over any base scheme). In contrast, the present paper focuses on computing the pure magnets (in the sense of the very recent framework algebraic magnetism) of the self-action of diagonalizable monoid schemes over an arbitrary base scheme. Recall that our goal in this paper is to show that algebraic magnetism captures sharpness and minimal generators modulo invertible elements, classical concepts in semigroup theory.

\ackn I thank the referee for the very helpful remarks. The author has received funding from the ISF grant No. 1577/23.
\xackn

 \section{Attractors and magnets for affine schemes}
 
 In this section, we recall some definitions and facts about Algebraic Magnetism for actions of diagonalizable monoids on affine schemes. 
 Let $S$ be a base scheme.
 Let $M$ be a cancellative abelian monoid and let $Z$ be $M^{gp}$. Let $N$ be a submonoid of $M$. The canonical inclusions of monoids $N \subset M \subset Z$ induce canonical morphisms of monoid schemes $A(Z)_S \to A(M)_S \to A(N)_S$. The diagonalizable group scheme $D(Z)_S$ identifies, as monoid scheme, with $A(Z)_S$. 
 
 Let $X\to S$ be a scheme endowed with an action of $A(M)_S$. The attractor $X^N$ is defined as the contravariant functor sending a scheme $T$ over $S$ to the set  $\Hom^{A(M)_T} ( A(N)_T , X_T )$. This is representable by a scheme in many cases, in particular if $X\to S$ is affine without any other assumption.  In this paper we will only deal with the case where $X$ is affine over $S$.  In fact we will only deal with the case where $X= A(L)_S$ for a submonoid $L \subset M$, in other words $a$ is the action of $A(L)_S$ induced by $A(M)_S \to A(L)_S$. Without loss of generality, we also assume that $S= \Spec (\mathbb{Z})$ (because attractors commute with base change and $a$ is defined over $\bbZ$ under the above assumptions). When $S = \mathbb{Z}$, we simply omit the subscript $\mathbb{Z}$ in notation. In particular $A(M)= \Spec (\mathbb{Z}[M])$.

\prop \label{propsv} If $X= A(L)$, then for any submonoid $N$ of $M$ the scheme $X^{N}$ is a closed subscheme of $X$ and equals $  \Spec (\mathbb{Z} [L] / J_N )$ where $J_N $ is the ideal of $\mathbb{Z} [L]$ given by $ \langle \mathbb{Z} X^m | m \in L \setminus (N\cap L) \rangle $ .
\xprop 
\pf
This is an immediate corollary of \cite[Theorem 3.20]{Ma}. In the following, for the convenience of the reader, we sketch the argument in this special case. It is enough to observe that for any ring $R$, we have canonical identifications \begin{align*}
X^N (R) &\cong \Hom ^{A(M)_R}(A(N)_R,X_R) \\ &\cong \Hom ^{A(M)_R} (A(N)_R,A(L)_R) \\ &\cong \Hom_{R\text{-alg}} ^{M\text{-graded}}(R[L], R[N]) \\
&\cong \Hom _{R\text{-alg} } ( R[L]/\langle R X^m | m \in L \setminus (N\cap L) \rangle, R)   \\
&\cong \Hom _{\bbZ\text{-alg} } ( \mathbb{Z}[L]/J_N, R)\\
& \cong \Spec (\mathbb{Z} [L] / J_N ) (R). 
\end{align*}
Let us justify the third identification  $\Hom_{R\text{-alg}} ^{M\text{-graded}}(R[L], R[N]) \cong \Hom _{R\text{-alg} } ( R[L]/\langle R X^m | m \in L \setminus (N\cap L) \rangle, R) $. Let $f \in \Hom_{R\text{-alg}} ^{M\text{-graded}}(R[L], R[N])$. Since it is a morphism of graded rings, $\langle R X^m | m \in L \setminus (N\cap L) \rangle$ is in the kernel of $f$. To obtain a morphism in $\Hom _{R\text{-alg} } ( R[L]/\langle R X^m | m \in L \setminus (N\cap L) \rangle, R)$ it is now enough to post-compose with the morphism $R[N] \to R, X^n \mapsto 1$. Conversely, let $F$ be a morphism in $\Hom _{R\text{-alg} } ( R[L]/\langle R X^m | m \in L \setminus (N\cap L) \rangle, R)$. We attach to $F$ the unique morphism in $\Hom_{R\text{-alg}} ^{M\text{-graded}}(R[L], R[N])$ sending, for all $m \in L $, $X^m \in R[L]_m$ to $F([X^m])X^m$ where $[X^m]$ is the class of $X^m$ in $R[L]/\langle R X^m | m \in L \setminus (N\cap L) \rangle$. This provides the desired reciprocal canonical bijections.
\xpf

The magnets of $a$ are by definition the submonoids of $M$. The set of magnets is denoted $m(a)$. A magnet $N$ is called pure if for any other magnet $P \subset N$ such that $X^P=X^N$, we have $P=N$. The set of pure magnets is denoted $\mho (a)$. If $N$ is a pure magnet, we put $m^N (a) = \{ N' \in m(a) | X^N = X^{N'} \}$.  We have a decomposition $m(a)= \sqcup _{N \in \mho (a) } m^N (a)$. If $N'$ is a magnet, we denote by $E(N')$ the pure magnet $N$ such that $N' \in m^{N}(a)$, i.e. the only pure magnet $N$ such that $X^{N'} = X^{N}$. Recall that $D(Z)_S$ acts also on $A(M)_S$ via $D(Z) \to A(M)_S$, let $a'$ be this action. The pure magnets of $a'$ are submonoids of $M$, and in fact we have  $\mho (a) = \mho (a')$, an equality of subsets of submonoids of $M$. The above sentence follows from Proposition \ref{propsv} or \cite[Proposition 3.30]{Ma}.
Recall that the goal of this paper is to compute $\mho (a)$, or equivalently $\mho (a')$.

\section{Pure Magnets of $A(N)_S$ where $N$ is cancellative and sharp}

Let $M$ be a cancellative (integral in Ogus) monoid (we do not assume that $M^{gp}$ is finitely generated). Let $E$ be a subset of $M$. Let $[E\rangle$ be the submonoid of $M$ generated by $E$. Let $a$ be the action of $A(M)_S$ on $X=A([E\rangle)_S$ corresponding to the inclusion $[E\rangle \subset M$.

\prop \label{PropaimantsAN}  Assume that \begin{enumerate} \item for all $E' \subsetneq E$, we have $[E'\rangle \subsetneq [E\rangle$ (i.e. $E'$ does not generate the monoid $[E\rangle$), and
\item $[E\rangle^*=0$,
\end{enumerate}
 then $\mho (a) \cong \{$subsets of $E\}$.
\xprop 
\pf We start with an observation.
\lemm \label{lemmBE} For all $B \subset E,$ for all $ B' \subsetneq B$,  we have $[B'\rangle \subsetneq [B\rangle$.
\xlemm
\pf
Let $B' \subsetneq B \subset E$ as in the statement. Assume by contradiction that $B'$ generates $[B\rangle$. Let $b \in B \setminus B'$. Then $b \in [B'\rangle$. This implies that $E \setminus \{b\}$ generates $[E\rangle$ which is a contradiction. 
\xpf
 We now prove Proposition \ref{PropaimantsAN}. We can assume that $S= \Spec (\mathbb{Z})$ by base change. Consider the map $f:\mho (a) \to \{$subsets of $E\}$, $N \mapsto N \cap E$. We claim that it is a  bijection. We start with injectivity. Let $N,N' \in \mho (a)$ such that $N \cap E = N' \cap E$.
Let us prove that $X^{N}=X^{[N \cap E\rangle}$. By Proposition \ref{propsv}, it is enough to prove that $J_{N}=J_{N\cap E}.$ Since $[N  \cap E \rangle \subset N$, we have $J_N \subset J_{[N \cap E \rangle }$. Let us prove conversely that $ J_{[N \cap E \rangle } \subset J_N $. It is enough to prove that $[E\rangle \setminus [N \cap E \rangle $ is included in the monoid-ideal of $[E\rangle$ generated by $[E\rangle \setminus N$. So let $e \in [E \rangle \setminus [N \cap E \rangle $. Put $V = N \cap E $ and $H = E \setminus (N \cap E)$, so that $E = V \sqcup H$. Since $e$ is outside $[V\rangle $, there exists $h \in H$ such that $e=e'+h$ with $e' \in [E\rangle$.   Since $H \cap N = \emptyset$, we observe that $H \subset [E\rangle \setminus N$. So $e=h+e'$ with $h \in [E\rangle \setminus N$. Therefore, $e$ belongs to the monoid-ideal of $[E\rangle$ generated by $[E\rangle \setminus N$. So the identity $X^{N}=X^{[N \cap E\rangle}$ is proved. We also have $X^{N'}=X^{[N' \cap E\rangle}$. So $X^N = X^{N'}$. Since $N $ and $N'$ are pure, this implies $N=N'$. This finishes injectivity. 
 To prove that $f$ is surjective, it is enough to prove that for all $B \subset E$, $[B\rangle$ is a pure magnet and $[B \rangle \cap E = B$. So let $B \subset E$ be a subset. Let us prove that $[B \rangle \cap E = B$. It is immediate that $B \subset [B\rangle \cap E \subset [B\rangle$. So $[B\rangle = [[B\rangle \cap E\rangle$. If, by contradiction, $B ~\subsetneq ~[B \rangle \cap E$, then by Lemma 3.2 $B$ does not generate  $[B\rangle = [[B\rangle \cap E\rangle$ which is absurd. So $[B \rangle \cap E = B$.  
It remains to prove that $[B\rangle$ is a pure magnet. 
So let $N \subset [B \rangle$ be a submonoid such that $X^N = X^{[B\rangle}$ and let us prove that $N = [B \rangle$. It is enough to prove that $B \subset N$. Assume by contradiction that there exists $b \in B $ such that $b \notin N$. Put $A= \mathbb{Z}[[E\rangle]$. Let $J_N \subset A$ be the ideal associated to $X^N$ and let $J_{[B\rangle}\subset A$ be the ideal associated to $X^{[B\rangle}$. Since $b \in [E\rangle \setminus N $, the element $X^b$ belongs to $J_N= \langle \bbZ X^e | e \in [E\rangle \setminus N \rangle$. Recall that $J_{[B\rangle}$ is the ideal $\langle \bbZ X^e | e \in [E\rangle \setminus [B \rangle \rangle $. The ideal $J_{[B\rangle}$ of $A$ corresponds to the monoid-ideal $\mathrm{MonIdeal}([E\rangle \setminus [B \rangle )$ of $[E\rangle$ generated by $[E\rangle \setminus [B \rangle $. Let us prove that $b \not \in \mathrm{MonIdeal}([E\rangle \setminus [B \rangle \rangle)$. Assume by contradiction that $b \in \mathrm{MonIdeal}([E\rangle \setminus [B \rangle \rangle)$. Then there exists $e' \in E \setminus B $, and for each $e \in E$ an integer $\lambda _e \in \mathbb{N}$ such that $b= \sum _{e \in E } \lambda _e e + e'$. 
If $\lambda _b=0$, then $b \in [E \setminus \{b\}\rangle$ which is absurd. If $\lambda _b \geq 1$ then since  $M$ is cancellative, we get the identity $0= \sum _{e \in E \setminus \{b \}} \lambda _e e + (\lambda _b -1) b +e'=0$. This implies that $e' \in N^*. $ 
So $N^* \ne 0$ which is absurd. We used that $E \cap \{0\}= \emptyset$ because $0$ has no effect on generation ($[H\rangle= [H \cup 0\rangle$ $\forall H$). So we proved that $b \not \in \mathrm{MonIdeal}([E\rangle \setminus [B \rangle \rangle)$.
So $X^b $ belongs to $J_N$ but does not belong to $J_{[B\rangle}$. So $J_N \ne J_{[B\rangle}$. So $X^N \ne X^{[B\rangle}$. This is a contradiction. So $B  \subset N$ and $N = [B \rangle$. So $[B\rangle$ is a pure magnet. This proves surjectivity and concludes the proof. 
\xpf 

Let us comment on the conditions of Proposition \ref{PropaimantsAN}. The first condition is of course necessary.
The second is also necessary, indeed let $M=\bbZ$ and $E = \{-1,1\}$. The set $E$ satisfies that for all $E' \subsetneq E$, $E'$ does not generate the monoid $[E\rangle$. The action $a$ is the translation action of $\mathbb{G}_m$ on itself. We have $\mho (a) = \{ 0, \mathbb{Z} \}$ which is of cardinality $2$ and not $4$.

\section{Pure magnets of $A(N)_S$ where $N$ is cancellative }

Let $M$ be a cancellative monoid. In the previous section, we explained how to compute the pure magnets of the action $a$ of $A(M)_S$ on itself if $M$ is sharp (i.e. if $M^*=0$). In this section, we provide results allowing to treat the general case from the sharp case. 

\lemm \label{group0M}Assume that $M$ is a group and $M \ne 0$. Then $\mho (a)= \{0, M \}$. Moreover $(A(M)_S)^N= \emptyset $ if $N \subsetneq M $ and $(A(M)_S)^M=A(M)_S$. If $M=0$, then $\mho (a)=\{0\}$.
\xlemm
\pf Let $N \subset M$ be a submonoid. If $N =M$, then $J_N =0$ so $(A(M)_S)^M=A(M)_S$. If $N \ne M$, then there exists $g \in M \setminus N$. Since $M$ is a group, $g$ has an inverse $-g$ in $N$. Then for all $n \in N $, we have $n = n+ (-g) +g  $ which shows that $J_N = \mathcal{O} _{A(M)_S}$, so $(A(M)_S)^N= \emptyset $. This implies that $\mho (a)= \{0, M \}$. The last assertion is trivial.
\xpf 

\prop \label{propjhg}Assume that $M^*\ne 0$. Let $N \subset M$ be a submonoid of $M$. If $M^* \cap N  \ne M^*$ (i.e. if $N$ does not contain $M^*$), then $(A(M)_S)^N = \emptyset. $
\xprop 
\pf  By \cite[Fact 3.15]{Ma} we have an equivariant morphism $A(M)_S \to A(M^*)_S$. We have $(A(M^*)_S)^N = A(M^*)_S^{M^* \cap N } $. So by Lemma \ref{group0M}, we have $(A(M^*)_S)^N= \emptyset$. This implies that $A(M)_S^N = \emptyset $ (the only scheme with underlying set zero is the empty scheme by \cite[\href{https://stacks.math.columbia.edu/tag/00E0}{Tag 00E0}]{stacks-project}).
\xpf 

Recall that the group $M^*$ acts on $M$ and that the quotient $M/M^*$ is a monoid. Let $\overline{a}$ be the action of $A(M/M^*)_S$ on itself. 
\prop \label{propNN'**} Let $N , N' $ be submonoids of $M$ containing $M^*$. Then 
$A(M)^N _S = A(M)^{N'}_S$ if and only if $A(M/M^*)_S^{N/M^*} = A(M/M^{*})_S ^{ N'/M^*}$.
\xprop 
\pf We reduce to the case where $S= \Spec( \mathbb{Z})$. Since $M^* \subset N$, we have a well-defined surjective map $M \setminus N \to M/{M^*} \setminus N /{M^*}$. Let $J_N $ and $ J_{N'} $ be the ideals of $\mathbb{Z}[M]$ corresponding to $A(M)^N _S $ and $ A(M)^{N'}_S$. Similarly,  let $J_{N/M^*} $ and $ J_{N'/M^*} $ be the ideals of $\mathbb{Z}[M/M^*]$ corresponding to $A(M/M^*)^{N/M^*} _S $ and $ A(M/M^*)^{N'/M^*}_S$. It is enough to prove that $J_N = J_{N'}$ if and only if $J_{N/M^*} = J_{N'/M^*} $. If $J_N = J_{N'}$, then $J_{N/M^*} = J_{N'/M^*} $. Reciprocally, assume that $J_{N/M^*} = J_{N'/M^*} $. It is enough to prove that $J_N \subset J_{N'}$. It is enough to prove that $M \setminus N$ is included in the monoid-ideal of $M$ generated by $M \setminus N'$.  So let $m \in M \setminus N $. Let $\overline{m}$ be its image in $M/M^* \setminus N / M^*$. By assumption, there exists $\overline{t} \in   M/M^* $ and $\overline{m'} \in M/M^* \setminus N' /M^*$ such that $\overline{m}=\overline{t}+ \overline{m'}$. This implies that there exists $t \in M , m' \in M\setminus N' $ and $g \in M^*$ such that $m= t+ g + m' $. So $m$ belongs to the monoid-ideal of $M$ generated by $M \setminus N'$. This finishes the proof. 
\xpf 
Let $f: M \to M/M^*$ be the projection. 
\lemm The maps $N \mapsto N/M^*$ and $P \mapsto f^{-1}(P)$ provides a bijection 
\[\{\text{submonoids of } N \text{ containing } M^* \} \leftrightarrow \{ \text{submonoids } P \text{ of } M/M^* \}  \]
\xlemm
\pf
The maps are well-defined and mutually inverse. 
\xpf  
\prop \label{propne0}Assume that $M^* \ne 0$, the set of pure magnets of the action $a$ of $A(M)_S$ on itself is $\mho (a)= \{0  \} \sqcup \{ f^{-1}(P)| P \in \mho (\overline{a})\}$.
\xprop 
\pf
Let $N $ be a submonoid of $M$. If $N$ does not contain $M^*$, then $A(M)_S^N=A(M)_S^0$ by Proposition \ref{propjhg}. This shows that the only pure magnet of $a$ not containing $M^*$ is $0$. Now let $N $ be a submonoid of $M$ containing $M^*$, we need to prove that $E(N)=f^{-1} (P)$ for some pure magnet $P$ of $\overline{a}$. Let $P$ be $E( N/M^*)$. In particular $A(M/M^*)_S^{P}= A(M/M^*)_S^{N/M^*}$. By Proposition \ref{propNN'**}, we have $A(M)_S^{f^{-1}(P)}=A(M)_S^{N}$. It is enough to prove that $f^{-1}(P)$ is a pure magnet of $a$. So let $Q \subset f^{-1} (P)$ be a submonoid such that $A(M)_S^Q= A(M)_S ^{f^{-1} (P)}$. We necessarily have $M^* \subset Q$ by Proposition \ref{propjhg}. By Proposition \ref{propNN'**} and because $P$ is a pure magnet, we have $f(Q)=P$. This implies $Q=f^{-1} (P)$. So $f^{-1} (P)$ is a pure magnet of $a$. 
\xpf

Recall that $M$ is an abelian and cancellative monoid.

\theo \label{finaltheo} Assume that $E$ is a minimal set of generators of $M/M^*$.
\begin{enumerate}
\item If $M^*=0$, then $(B \subset E) \mapsto [B\rangle $ provides a bijection between  $\mathscr{P} (E)$ and $\mho (a)$. 

\item If $M^*\ne 0$, then $\mho (a) =\{0\} \sqcup \{ f^{-1} ([B\rangle)| B \in \mathscr{P}(E)\}$, where $f: M \to M/M^*$. 
\end{enumerate} 
In particular $\#\mho(a)= 2^{\#E}$ if $M^*=0$ and $\#\mho(a)= 2^{\#E}+1 $ if $M^*\ne 0$.

\xtheo 
\pf If $M^*= 0$ this is Proposition \ref{PropaimantsAN}.
Otherwise, \cite[Proposition 1.3.3]{Og18} shows that the sharp monoid $M/M^*$ is cancellative (integral). Therefore it is enough to apply Propositions \ref{PropaimantsAN} and \ref{propne0}. The last assertion is immediate.
\xpf 
\section{Examples}

\begin{enumerate}
\item Let $a$ be the action of $D(\bbZ) $ on $A(\bbN)$ corresponding to $\bbN \subset \bbZ$, geometrically this is the multiplication of $\mathbb{G}_m$ on $\bbA^1$. Then since $\bbN $ is sharp and $\{1\}$ is the minimal set generating $\mathbb{N}$, Theorem \ref{finaltheo} implies that $\mho (a) = \{ 0 , \bbN \}$. 
\item Let $a$ be the action of $D(\bbZ \times \bbZ) $ on $A(\bbZ \times \bbN) $, i.e. geometrically $\bbG_m ^2 $ acting on $\bbG_m \times \bbA^1$ by multiplication. Then $(\bbZ \times \bbN)^*=\bbZ \times 0$ and $(\bbZ \times \bbN)/(\bbZ \times \bbN)^*\simeq \bbN$, so Theorem \ref{finaltheo} implies that $\mho (a)= \{ 0, \bbZ , \bbZ \times \bbN \}$.  Note that $0 \times \bbN $ is not a pure magnet. We have $A(\bbZ \times \bbN)^{0 \times \bbN}= \emptyset$ (empty scheme) because $0 \times \bbN$ does not contain $(\bbZ \times \bbN)^*$.
\item Let $a$ be the action of $D( \mathbb{Z}^2) $ on $A ([(1,-1),(1,0),(1,1)\rangle$. Then since $(1,-1),(1,0)$ and $(1,1)$ are independent, and $[(1,-1),(1,0),(1,1)\rangle$ is sharp, Theorem \ref{finaltheo} implies that $\# \mho (a) = 2^3 $. 
 
\item Let $a$ be the action of $D(\bbZ \times \bbZ /2\bbZ )$ on $A(N)$ where $N= [(0,\overline{1}) , (2, 0) , (3, \overline{1})  \rangle$. Then $N^* = \{(0,0),(0, \overline{1})\} \cong \bbZ/2\bbZ$. We have $N /N^* \simeq [2,3\rangle \subset \bbN $. The set $\{2,3\}$ is a minimal generating set of $N /N^* \simeq [2,3\rangle$. So $\mho (a) = \{ 0, [(0, \overline{1})\rangle, [(0, \overline{1}), (2,0)\rangle , [(0, \overline{1}), (3,0)\rangle, N \rangle \}$.

\item Let $a$ be the action of $D(\mathbb{Z}^2)$ on $A([\{ (k,2k-1) |k \in \mathbb{N} \}\rangle)$. Observe that $\{ (k,2k-1) |k \in \mathbb{N} \}$ is a minimal set of generators of $[\{ (k,2k-1) |k \in \mathbb{N} \}\rangle$ (which is sharp). Indeed, if $n_1+\ldots +n_s= n$ with $s>2$ (equality of integers in $\mathbb{N} _{\geq 1 }$), then $((2n_1-1+\ldots +2n_s-1)<2 n-1$.
Therefore Theorem \ref{finaltheo} implies that $\mho (a) $ is in bijection with $\mathscr{P} (\{ (k,2k-1) |k \in \mathbb{N} \})$, in particular $\#\mho (a) >\aleph _0.$
\item Let $N $ be a submonoid of $\mathbb{N}$, and let $a$ be the action of $ A(\bbN)$ on $A(N)$. Then $\mho (a)$ is finite by Theorem \ref{finaltheo}. Indeed any numerical monoid (a submonoid of $\bbN$) is finitely generated. For example, if $N=d \mathbb{N}$, then $N^*=0$ and $\{d\}$ is a minimal set of generators of $N/N^*\cong N$ so that the only pure magnets are just $0$ and $N$ by Theorem \ref{finaltheo}. In the following, we provide a direct proof that the pure magnets are just $0$ and $N$ in this case.  Assume that $0 \subsetneq P \subsetneq d \mathbb{N}$ is a strict submonoid, then $d$ is outside $P$. So $\Hom ^{N\text{-graded}}_{R\text{-alg}}(R[N],R[P])$ is a singleton for all rings $R$. So $X^P=X^0$ and $P$ is not a pure magnet. It is now established that $0$ and $N$ are the only pure magnets because $A(N)^N = A(N)\ne A(N)^0$.
\end{enumerate}

\end{document}